\input amstex

\documentstyle {amsppt}
\magnification=\magstep1 \NoRunningHeads

\topmatter
\title
A Note on the Zeros of Jensen Polynomials
\endtitle
\author
Young-One Kim and Jungseob Lee
\endauthor

\abstract A recent result of Griffin, Ono, Rolen and Zagier on
Jensen polynomials related with the Riemann zeta function is
improved.
\endabstract

\address
Department of Mathematical Sciences and Research Institute of
Mathematics, Seoul National University, Seoul  08826,  Korea
\endaddress

\email kimyone\@snu.ac.kr
\endemail

\address
Department of Mathematics, Ajou University, Suwon 16499,  Korea
\endaddress

\email jslee\@ajou.ac.kr
\endemail



\thanks
The first author was supported by the SNU Mathematical Sciences
Division for Creative Human Resources Development.
\endthanks

\date 2021. 5. 12.  \enddate

\endtopmatter

\document

This note is concerned with the zeros of entire functions. Let $f$
be an entire function. We denote the zero set of $f$ by $\Cal Z(f)$,
that is, $\Cal Z(f)=\{z\in\Bbb C: f(z)=0\}$; the order
$\rho=\rho(f)$ of $f$ is defined by
$$
\rho=\limsup_{r\to\infty}\frac{\log\log M(r;f)}{\log r},
$$
where $M(r;f)=\max_{|z|=r}|f(z)|$; and $f$ is said to be {\it real}
if $f^{(n)}(0)\in\Bbb R$ for all $n$.

Throughout this note $f$ denotes a real entire function of order
$\rho<2$, and in the case where $f$ is even $f_0$ denotes the
function such that $f(z)=f_0(z^2)$ for all $z$. In this case, $f_0$
is a real entire function of order $\rho/2<1$, $\Cal Z(f_0)=\{ z^2:
z\in\Cal Z(f)\}$, and we have
$$
f_0^{(k)}(0)=\frac{k!}{(2k)!}f^{(2k)}(0) \qquad (k=0, 1, 2, \dots).
$$
Since $f$ is a real entire function, we have $f(z)=0$ if and only if
$f(\bar z)=0$; and since $f$ is of order less than $2$, it can be
represented in the form
$$
f(z)=cz^me^{bz}\prod_j\left(1-\frac{z}{a_j}\right)e^{z/a_j},
$$
where $c,b\in \Bbb R$, $m$ is a non-negative integer and $\sum
|a_j|^{-2}<\infty$. As a consequence, there are real polynomials
$P_1, P_2, \dots$ such that $\Cal Z(P_k)\subset\Cal Z(f)\cup\Bbb R$
for all $k$, and $P_k\to f$ uniformly on compact sets in the complex
plane. In particular, if $\Cal Z(f)\subset\Bbb R$ then $\Cal
Z(f^{(n)})\subset\Bbb R$ for every $n$ unless $f^{(n)}$ is
identically equal to zero. The Riemann Xi-function is given by
$$
\Xi(iz)=\frac{1}{2}\left(z^2-\frac{1}{4}\right)\pi^{-\frac{z}{2}-\frac{1}{4}}
\Gamma\left(\frac{z}{2}+\frac{1}{4}\right)\zeta\left(\frac{1}{2}+z\right).
$$
It is an even real entire function of order $1$, its zero set is
contained in the infinite strip $\Bbb S=\{z\in\Bbb C:
|\operatorname{Im} z|<1/2\}$, and the Riemann hypothesis is  the
conjecture that $\Cal Z(\Xi)\subset\Bbb R$, which is equivalent to
that $\Cal Z(\Xi_0)\subset [0, \infty)$.

A real polynomial is said to be {\it hyperbolic} if it has real
zeros only. The Jensen polynomials $J(f;1), J(f;2), \dots $ of $f$
are defined by
$$
J(f;d)(z)=\sum_{k=0}^d\binom dk f^{(k)}(0)z^k\qquad (d=1, 2, \dots).
$$
It is known \cite{3,8} that $\Cal Z(f)\subset\Bbb R$ if and only if
all the Jensen polynomials of $f$ are hyperbolic. Let $d$ be a
positive integer. If there is a non-negative integer $N$ such that
$J(f^{(n)};d)$ is hyperbolic for every $n\geq N$, we denote the
smallest  such integer by $N(f;d)$, otherwise we write
$N(f;d)=\infty$. One can easily check that the Riemann hypothesis
holds if and only if $N(\Xi_0; d) = 0$ for all $d$. Griffin et al.
\cite{2} showed that $N(\Xi_0;d)<\infty$ for all $d$ by establishing
an asymptotic formula for the sequence $\langle J(\Xi_0^{(n)};
d)\rangle_{n=1}^{\infty}$ which relates (for large $n$) the zeros of
$ J(\Xi_0^{(n)}; d) $ with those of the $d$-th Hermite polynomial.
We improve the result by applying some general theorems on real
entire functions of order less than $2$ with restricted zeros.
However it should be remarked that our proof gives no information
about the Jensen polynomials considered, except that they are
hyperbolic.

\proclaim{Theorem 1} Suppose that $f$ is a transcendental real
entire function of order $\rho<2$ and $\Cal Z(f)\subset\Bbb S$. Then
for every $c>\rho$ we have $N(f;d)=O(d^{c/2})$ as $d\to\infty$. If
$f$ is even, we also have $N(f_0;d)=O(d^{c/2})$ as $d\to\infty$ for
every $c>\rho$.
\endproclaim

Since $\Xi$ is of order $1$, it follows that $N(\Xi;d)=O(d^c)$ and
$N(\Xi_0;d)=O(d^c)$ as $d\to\infty$ for every $c>1/2$.   Theorem 1
is a consequence of Theorems 2 and 3 below. Let $ \tilde\Bbb
S=\{z^2: z\in\Bbb S\}=\{z\in\Bbb C: (\operatorname{Im} z)^2 -
4^{-1}< \operatorname{Re} z \}$.

\proclaim{Theorem 2}  Suppose that $f$ is a transcendental real
entire function of order $\rho<2$ and $\Cal Z(f)\subset\Bbb S$. Then
for every $c>\rho$ there is a positive integer $n_1$ such that
$$
\Cal Z(f^{(n)})\subset \left\{z\in\Bbb S: |\operatorname{Re} z|\geq
n^{1/c}\right\}\cup\Bbb R \qquad (n\geq n_1). \tag 1
$$
If $f$ is even, then for every $c>\rho$ there is a positive integer
$n_2$ such that
$$
\Cal Z(f_0^{(n)})\subset \left\{z\in\tilde\Bbb S: \operatorname{Re}
z\geq n^{2/c}\right\}\cup [0, \infty) \qquad (n\geq n_2).
$$
\endproclaim

\demo{Proof} The first part follows from Theorem 2 in \cite{4}, and
the second part is a consequence of Theorem 1 in \cite{5} (with
$k=1/2$) applied to $g$ given by $g(z)=f_0(z-4^{-1})$. \qed
\enddemo

For $\delta>0$ we set $S(\delta)=\{z\in\Bbb C: |\operatorname{Im}
z|\leq \delta |z|\}$.

\proclaim{Theorem 3} Suppose that $P$ and $Q$ are real polynomials,
$\delta>0$, $\Cal Z(P)\subset S(\delta)$, $Q$ is hyperbolic and
$\deg Q\leq \delta^{-2}$. Then the polynomial
$$
P(D)Q=\sum_{k=0}^{\deg P}\frac{P^{(k)}(0)}{k!} Q^{(k)}
$$
is hyperbolic.
\endproclaim

\demo{Proof} See \cite{6}. \qed
\enddemo

\proclaim{Corollary}  Suppose that $P$ is a  real polynomial,
$\delta>0$ and $\Cal Z(P)\subset S(\delta)$.  Then $J(P;d)$ is
hyperbolic for $d\leq \delta^{-2}$.
\endproclaim

\demo{Proof} If $Q(z)=z^d$, then $J(P;d)(z)=z^d (P(D)Q)(z^{-1})$.
\qed
\enddemo

\demo{Proof of Theorem 1} Suppose that $c>\rho$. By Theorem 2, there
is a positive integer $n_1$ such that (1) holds. Let $d$ and $n$ be
positive integers such that
$$
n\geq\max\{n_1, (d/4)^{c/2}\},
$$
and put  $\delta=2^{-1}n^{-1/c}$. Then $d\leq \delta^{-2}$, and from
(1) we can see that $|\operatorname{Im} z|\leq \delta |z|$ for all $
z\in \Cal Z(f^{(n)})$. Hence we have $\Cal Z(f^{(n)})\cup\Bbb
R\subset S(\delta)$. Let $P_1, P_2, \dots$ be real polynomials such
that $\Cal Z(P_k)\subset\Cal Z(f^{(n)})\cup\Bbb R$ for all $k$, and
$P_k\to f^{(n)}$ uniformly on compact sets in the complex plane.
Then $J(P_k;d)\to J(f^{(n)};d)$ uniformly on compact sets in the
complex plane, and the corollary to Theorem 3 implies that
$J(P_k;d)$ is hyperbolic  for every $k$. Therefore $ J(f^{(n)};d)$
is hyperbolic, and we conclude that
$$
N(f;d)\leq \left\lceil\max\{n_1, (d/4)^{c/2}\}\right\rceil \qquad
(d=1, 2, \dots).
$$
This proves the first part of the theorem. The second part is
analogously proved. \qed
\enddemo

\remark{Remark} Since our argument deals with entire functions of
order less than $2$, it refines the proof of Theorem 3.6 in
\cite{1}.
\endremark

Suppose that $T>0$ and  $\Cal Z(f)\subset\Bbb S(T)= \{z\in\Bbb S :
|\operatorname{Re} z|\geq T\}\cup\Bbb R$. Then the same argument as
in the proof of Theorem 1 shows that $J(f;d)$ is hyperbolic whenever
$d\leq 1+4T^2$, because we have $\Bbb S(T)\subset S(\delta)$ with
$\delta=(1+4T^2)^{-1/2}$. Since the condition  $\Cal Z(f)\subset\Bbb
S(T)$ does not imply $\Cal Z(f')\subset\Bbb S(T)$ in general, it is
not clear whether the polynomials $J(f^{(n)}; d)$, $n=1, 2, \dots$,
are hyperbolic for $d\leq 1+4T^2$. If $f$ is even and $T\geq 1/2$,
then it is the case.

\proclaim{Theorem 4} Suppose that $f$ is an even transcendental real
entire  function of order less than $2$, $T\geq 1/2$ and $\Cal
Z(f)\subset \Bbb S(T)$. Then
 $J(f^{(n)}; d)$ is
hyperbolic for $d\leq 1+4T^2$ and for every $n$.
\endproclaim

\demo{Proof} Let $\delta=(1+4T^2)^{-1/2}$,
$\tilde\delta=2\delta\sqrt{1-\delta^2}\,$($=T^{-1}(1+4^{-1}T^{-2})^{-1}$)
and $\tilde S=\{z\in S(\tilde\delta): \operatorname{Re} z\geq 0\}$.
Then $\Bbb S(T)\subset S(\delta)$, and $\tilde S$ is convex. Since
$T\geq 1/2$, we have $0<\delta\leq 2^{-1/2}$, and it follows that
$z\in S(\delta)$ if and only if $z^2\in\tilde S$.

We prove the theorem by showing that $\Cal Z(f^{(n)})\subset
S(\delta)$ for every $n$. We have $\Cal Z(f)\subset\Bbb S(T)\subset
S(\delta)$, and this implies that $\Cal Z(f_0)\subset \tilde S$.
Since $f_0$ is of order less than $1$ and $\tilde S$ is convex, the
Gauss-Lucas theorem implies that $\Cal Z(f_0')\subset\tilde S$.
Hence $\Cal Z(f')\subset S(\delta)$, because $f'(z)=2zf_0'(z^2)$. We
have
$$
f''(z)=2\left(f_0'(z^2) + 2z^2 f_0''(z^2)\right).
$$
If we write $g(z)=z(f_0'(z))^2$, then $g$ is of order less than $1$
and $\Cal Z(g)\subset \tilde S$. Again, it follows from the
Gauss-Lucas theorem that $\Cal Z(g')\subset \tilde S$, that is, all
the roots of the equation
$$
f_0'(z)\left(f_0'(z) + 2z f_0''(z)\right)=0
$$
lie in $\tilde S$. Therefore $\Cal Z(f'')\subset S(\delta)$. Since
$f''$ is even,   the result follows from an induction. \qed
\enddemo

\remark{Remark} The proof shows that $\Cal Z(f_0^{(n)})\subset\tilde
S\subset S(\tilde\delta)$ for all $n$. Hence $J(f_0^{(n)};d)$ is
hyperbolic for $d\leq T^{2}(1+4^{-1}T^{-2})^{2}$ and for every $n$.
\endremark

Currently, it is known \cite{7} that $\Cal Z(\Xi)\subset \Bbb S(T)$
for some $T\geq 3\cdot 10^{12}$. Thus $J(\Xi^{(n)};d)$ is hyperbolic
for $d\leq 3.6\cdot 10^{25}$ and for every $n$, and
$J(\Xi_0^{(n)};d)$ is hyperbolic for $d\leq 9\cdot 10^{24}$ and for
every $n$.


\Refs

\widestnumber\key{4}



\ref \key 1\by  M. Chasse \pages 875--885\paper Laguerre multiplier
sequences and sector properties of entire functions\yr 2013 \vol
58\jour Complex Var. Elliptic Equ.
\endref

\ref \key 2\by M. Griffin, K. Ono, L. Rolen and  D. Zagier \pages
11103--11110\paper Jensen polynomials for the Riemann zeta function
and other sequences\yr 2019 \vol  116 \jour  Proc. Natl. Acad. Sci.
USA
\endref



\ref \key 3\by M.-H. Kim, Y.-O. Kim and J. Lee  \pages
1121--1132\paper
 The convergence of a
sequence of polynomials with restricted zeros\yr 2019 \vol 478 \jour
J. Math. Anal. Appl.
\endref

\ref \key 4\by Y.-O. Kim \pages 819--830\paper Critical points of
real entire functions and a conjecture of Polya\yr 1996 \vol
124\jour  Proc. Amer. Math. Soc.
\endref

\ref \key 5\bysame \pages 657--667\paper Critical zeros and nonreal
zeros of successive derivatives of real entire functions\yr 1996
\vol 33\jour J. Korean Math. Soc.
\endref



\ref \key 6 \by  N.   Obreschkoff\pages  175--184 \paper Sur une
g\'{e}n\'eralisation du th\'eor\`eme de Poulain et Hermite pour les
z\'eros r\'eels des polyn\^omes r\'eels \yr 1961 \vol  12 \jour Acta
Math. Acad. Sci. Hungar.
\endref

\ref \key 7 \by  D. Platt and T. Trudgian  \paper The Riemann
hypothesis is true up to $3\cdot 10^{12}$   \jour Bull. London Math.
Soc.  \toappear
\endref


\ref \key 8\by G. P\'olya and I. Schur \pages 89--113 \paper
\"{U}ber zwei Arten von Faktorenfolge in der Theorie der
algebraischen Gleichungen \yr 1914 \vol 144\jour J. Reine Angew.
Math.\endref



\endRefs
\enddocument